\documentclass[12pt,reqno]{article}
\usepackage{amsmath,amssymb,amsthm}
\usepackage{comment}
\usepackage{float}
\usepackage{hyperref}
\usepackage[russian,english]{babel}
  \usepackage[utf8]{inputenc}
\usepackage{graphicx}
\usepackage[font=small]{caption}
\usepackage{setspace}
\usepackage{url}

\newcommand{\Qn}{\medskip\noindent\textbf{Q:\ }}
\newcommand{\Ar}{\medskip\noindent\textbf{A:\ }}

\newcommand\eps{\ensuremath{\varepsilon}}
\newcommand\Fraisse{Fra\"\i ss\'e}
\renewcommand\phi{\ensuremath{\varphi}}
\newcommand\R{\ensuremath{\mathbb R}}
\newcommand\ru{\foreignlanguage{russian}}
\newcommand\s{\mathcal S}

\title{Novosibirsk algebra and logic in the mid 1960s:\\
 A personal perspective}
\author{Yuri Gurevich\\
\small Computer Science and Engineering,
\small University of Michigan, USA}
\begin{document}
\date{}
\maketitle
\thispagestyle{empty}

\begin{center}
To the memory of Asan Dabsovich Taimanov
\end{center}

\bigskip
\begin{quote}\raggedleft\small\it
Of all that we experience,\\
there eventually of course remains only a memory,\\
but just in this way all lasting things
retain some of their actuality.\\[1.ex]
--- Kurt G\"odel \cite[\S1.3]{Wang}
\end{quote}

\begin{figure}[H]
\begin{center}
\includegraphics[width=9cm]{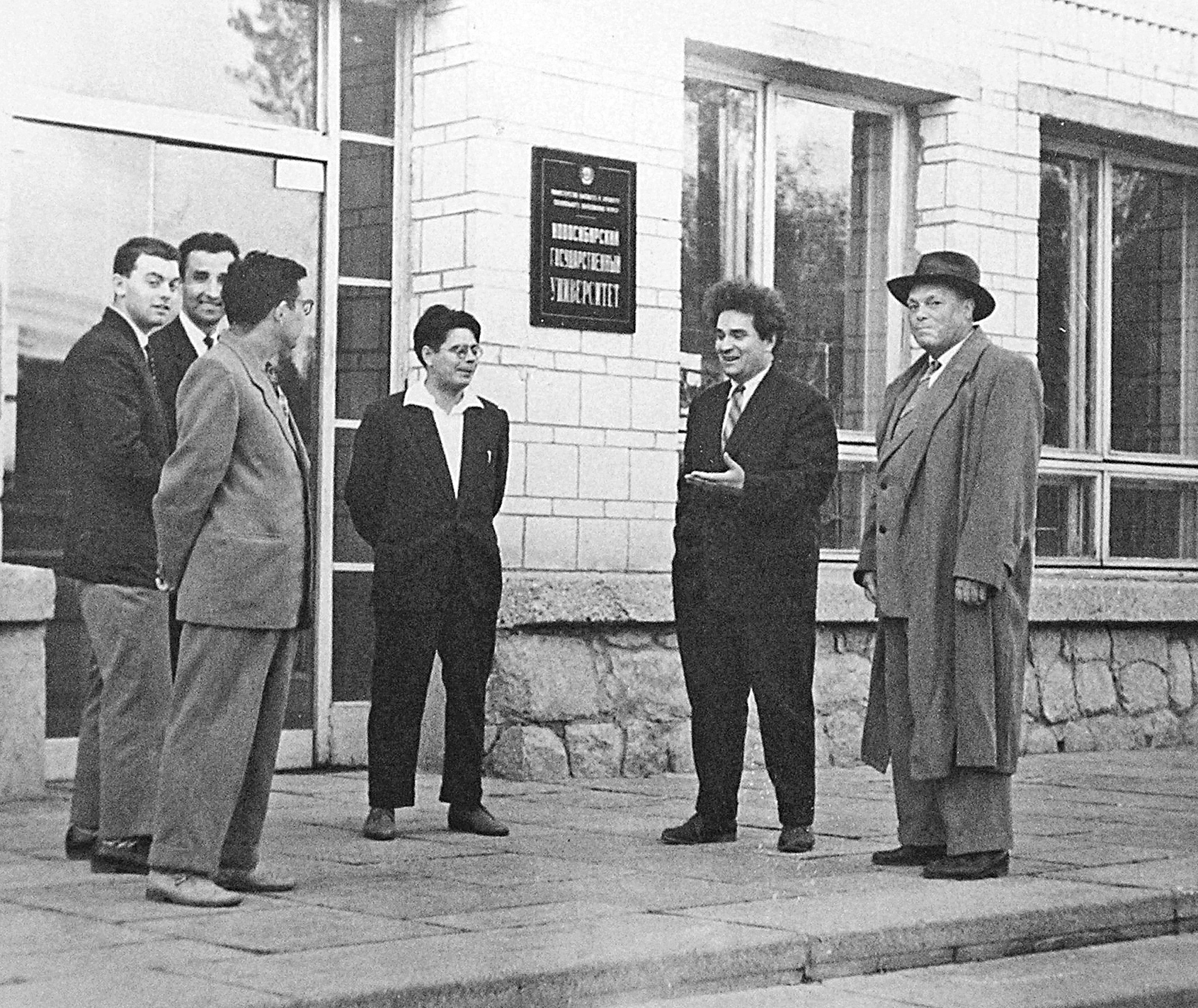}
\end{center}
\caption{\small The author, Ghaybulla Salikhov, Asan Taimanov, Mikhail Kargapolov,  Pavel Belinsky, and Anatoly Malcev}
\label{Pic}
\end{figure}

\noindent
\begin{quote}
Being asked to write about Asan D.\ Taimanov, I had a little problem. I knew Taimanov. I liked and respected him, and so I wanted to write something. But I didn't know him well. I didn't live in Novosibirsk\footnotemark where he lived and worked; I just visited Novosibirsk a number of times.
So I decided to widen the topic and write about my own Novosibirsk experience, including interactions with Taimanov.
\end{quote}

\footnotetext{Here and elsewhere in this article, speaking about Novosibirsk, I mean the so-called Akademgorodok \cite{W-ag}, the home of Novosibirsk State University and the Siberian Branch of the Russian Academy of Sciences.}

\section{My interactions with Taimanov}

\bigskip\noindent\textbf{Quisani%
\footnote{A former student of the author at the University of Michigan.}: }%
How well did you know Taimanov?

\Ar For me, Taimanov was a talented mathematician who independently discovered the so-called back-and-forth method in model theory \cite{Taimanov}, coauthored the influential article \cite{ELTT} on elementary theories, etc.
I respected him very much.
But all our conversations were mathematical.
We never had a serious conversation on anything else.

\Qn Why not?

\Ar We were of different generations, and I did not spend that much time in Novosibirsk where he was.
I remember several Kazakh students of Taimanov; he was Kazakh as well.
The students were my age, and I spoke more to them.

\Qn Was Taimanov approachable?

\Ar Yes. He was friendly, polite, and rather modest.
As far as I know, he was a good man, which is a nontrivial thing to say about the Novosibirsk mathematical scene of the second part of the previous century.

\Qn Whoa, this is serious.
We'll get there eventually I am sure.
For the time being, let me stick to lighter stuff.
Do you remember any funny story involving Taimanov?

\Ar I remember one incident at the 1983 Logic, Methodology and Philosophy of Science (LMPS) congress in Salzburg, Austria.
Most Soviet participants were afraid to talk to me.
Having left the land of victorious socialism in 1973, I had become a traitor of sorts.
But Taimanov and I had a little mathematical conversation.
As I was explaining some counterintuitive result,
the English phrase ``Do you buy that?'' was on my mind, and
I translated it to Russian literally.
Taimanov became visibly uncomfortable.
Apparently, it sounded to him as if I proposed some financial transaction.
I hurried to clarify that all I meant was ``Do you believe that?'' It was a bit awkward.

\Qn It wasn't the first time you met former Soviet colleagues?

\Ar Definitely not.
I remember the previous, 1979 LMPS congress, in Hannover, Germany.
There, in the last night, there was a party with free drinks, but food wasn't free.
The Soviet participants didn't have much Western currency to begin with, and they had virtually none by the last night.
Many of them became drunk, and some were willing and even eager to talk to me.
One of them, let me call him A, told me confidentially that another one, B, is the main KGB representative in their group.
Later B told me confidentially that A is a KGB man.

\Qn It is hard for me to understand Soviet life.
Taimanov was Kazakh, and you were Jewish. Did that create a distance between you in Russia?

\Ar No. If anything, that was bringing us a bit closer.
Imperial Russian chauvinism grew stronger and stronger after the Second World War.
In Russia, both Kazakhs and Jews were (marginalized as) national minorities.

\Qn Could you discuss national issues with your Russian friends?

\Ar Sure.
Russia has traditionally a stratum of open minded and often dissident intelligentsia.
My friends typically belonged to that wonderful stratum.
At the time, there were Jews and those with just one Jewish  parent in Russia who tried to hide their Jewishness.
I was a Jewish patriot. I believe this made it easier for me to discuss national issues with my Russian friends.

\section{Provincial in an academic center}

\Qn What was your connection to Novosibirsk University?
I understand that you didn't study there.

\Ar No, I didn't. I grew up in the provincial city of Chelyabinsk and, in 1957, entered the local Polytechnic.
Mathematical education there could have been better.
A couple of times I attempted to say to my math professor that the alleged proof of his wasn't really a proof, and each time he made everybody laugh at me.
In 1959, I succeeded to transfer to the math division of Ural State University, known by its Russian acronym UrGU.
The most active seminar there was in group theory, and so I chose to specialize in group theory.

\Qn How did you like it?

\Ar In the beginning I loved it, but then I grew somewhat disillusioned. It was abstract group theory, virtually divorced from classical mathematics and from applications. It was relatively easy to formulate problems of unquestionable difficulty but questionable importance.

\Qn How did you become interested in logic?

\Ar I told this story to Cris Calude \cite[\#213]{Gur}.
Here's a shorter version.
In late 1962, when I was a first-year graduate student, Mikhail I.\ Kargapolov, of Novosibirsk State University (NSU) gave a talk at UrGU.
He presented his simplification \cite{Kargapolov1} of  Wanda Szmielew's classification of abelian groups by elementary properties \cite{Szmielew}.

\Qn What does it mean to classify by elementary properties?

\Ar Consider a structure $X$, e.g.\ an abelian group. In second-order logic, you can quantify (variables ranging) over sets of elements of $X$, but in first-order logic, you can quantify only over elements of $X$.
Accordingly, first-order properties of $X$ are called elementary.

To classify abelian groups by elementary properties means to assign numerical parameters to the abelian groups in such a way that two groups are elementarily equivalent (have the same elementary properties) if and only if they have the same values of the parameters.

\Qn Define a $\{0,1\}$ parameter for every elementary property, and there you have a classification.

\Ar This is a good point. The desired parameters should be natural.

\Qn Naturality is in the eye of the beholder.

\Ar As an algebraist, I felt that Szmielew's parameters were not only natural, but the right ones.
Besides, her classification led to a decision algorithm for the elementary theory of abelian groups: Given an elementary property \phi, the algorithm decides whether \phi\ holds in all abelian groups or it fails in some.
At the time though, I had not heard of decision algorithms, and Kargapolov did not mention them.

An older student, Ali Kokorin, worked on ordered groups, and he
asked Kargapolov about ordered abelian groups (OAGs).
Kargapolov dodged the question, and Kokorin whispered to me: ``He is working on that.
Let's trounce%
\footnote{This is a polite translation.} him.''

We proved that no universal first-order formula distinguishes between any nontrivial (different from 1) OAGs and sent the result%
\footnote{eventually published in Algebra and Logic \cite{G002}}
to the national algebra conference scheduled for spring 1963 in Novosibirsk.
Then I continued the investigation alone; Kokorin was finishing his dissertation.
Little by little I kept moving forward, deepening my understanding of OAGs. By the time of the conference, I classified OAGs by elementary properties.
The parameters seemed right to me, but the proof was rather involved.

\Qn Could you present the classification result at the conference?

\Ar Yes. Kargapolov gave an hour-long talk.
He classified finite-rank OAGs.
Anatoly I.\ Malcev, a member of the Soviet Academy of Sciences, praised his achievement.
Right after that, I had 10 minutes for my talk, presumably to present the Gurevich--Kokorin result.
Instead, I spoke on my classification of OAGs.
Nobody seemed to listen; there were a lot of people and a lot of noise.
Suddenly, somebody\footnote{Mikhail Taitslin, as I learned later} shouted: ``So Kargapolov's result is a trivial special case of yours?''
It became quiet, and I chose my words carefully.
``No, Kargapolov's result is important, but yes, it is a special case.''

Malcev asked me to explain my proof at his seminar. It did not go well.
As I tried to explain things carefully, they were impatient.
I did not know model theory.
The little model theory needed for my proof, I developed myself, in the narrow frame of OAGs.
For me, it was a part of my work. But it wasn't new to them.
Then I would accelerate. But they weren't experts in OAGs, and I  would lose them.
At the end, they were skeptical, and Malcev advised me to write things down carefully. ``If the proof survives, come back.''

I came back in early fall.
During the summer, I wrote down my classification.
Also, assuming a decision algorithm for the elementary theory of linearly ordered sets \cite{Ehrenfeucht1,LL}, I constructed a decision algorithm for the elementary theory of OAGs.
Malcev asked Yuri Ershov to check my proofs.
For several days, I was explaining my stuff to Ershov.
Eventually he was convinced that it all works.

Malcev invited the two of us plus Kargapolov and Taimanov for a dinner at his house. The atmosphere was stilted. Malcev did almost all of the talking, mostly about famous people. He told us for example that the Japanese emperor had a doctoral degree in biology.
Feeling somewhat provincial, I asked Malcev what kind of meat we were eating. ``Yak,'' he said and continued talking about higher things. Feeling that I had already exposed my provinciality, I continued:``How do you know that it is yak?''``It came from the distribution outlet for academicians,'' he replied, ``and it is hard to chew. So what can it be?''

\Qn A distribution outlet for academicians? I don't understand.

\Ar Soviet society was hierarchical, and nowhere was the hierarchy more obvious than in Novosibirsk \cite{W-ag}.
The academicians (full members of the Academy of Sciences) were at the top of the academic hierarchy and had exclusive access to a special distribution outlet.

\Qn What happened next with your OAG results?

\Ar A few months later, I defended, in Novosibirsk, my Candidate dissertation%
\footnote{See \cite{W-CofS}}, based on my OAG work \cite[\#3]{Gur}. In a sense, I belonged to Malcev's scientific school, and he supported the defense. In fact he chaired the dissertation committee.

\Qn Did you transfer to Novosibirsk University?

\Ar No. According to the Soviet academic system, you defended your thesis at an academic center with sufficient expertise in the relevant area; it did not have to be your home institution.

\Qn Interesting. I like that.

\Ar Actually, Malcev offered to become formally my adviser, but I didn't want to offend Petr G.\ Kontorovich, my UrGU adviser. ``I'd be happy to join your group after the defense,'' I said.

Later Misha (Mikhail) Taitslin told me that declining the offer to become Malcev's student was a foolish mistake and that I was the first Jew whom Malcev was willing to advise.
``What about you?'' I asked.
Taitslin was Jewish and worked in Malcev's group.
``I imposed myself on Malcev," he said.

Taitslin predicted that I would be unable to join Malcev's group, and he turned out to be right.
Malcev kept telling me that he wants to hire me but can't do it right now.
I knew that it was mere politeness; Malcev was super polite in general.
Once, I decided to validate my understanding.
I asked a human resources ({\ru{отдел кадров}) administrator what would it take to join Malcev's group.
``Just bring me his signature,'' she said, ``on any piece of paper.''

\Qn So Malcev lied to you.

\Ar Lied is a strong word. I think that he could hire me if he wanted, but he might have his constraints.
The Moscow head of Steklov Mathematical Institute, Ivan M. Vinogradov, was openly antisemitic.
Antisemitism was growing in Novosibirsk as well.
Every Jewish dissertation defense in Novosibirsk had negative votes.
Malcev openly supported my Candidate defense, and all the reviews were strongly positive.
Yet, four members of my dissertation committee (out of about a dozen) voted negatively.
``This,'' Malcev said to me, ``is not against you but against me.'' That was undoubtedly a white lie.

\Qn Have you spent much time in Novosibirsk?

\Ar No. To arrange the thesis defense, I spent a few weeks in Novosibirsk.
After that I was coming there mostly to speak at the Algebra and Logic Seminar, but each time I tried to spend there a few days.

Novosibirsk was a great mathematical center.
Their math library was superb.
I mentioned Malcev already. He was an outstanding mathematician.
I remember one seminar when Malcev couldn't quite understand the speaker and kept asking him questions.
That by itself was not unusual. Malcev wasn't the quickest understander.
But this time his questioning went for much longer than usual.
Everybody else seemed to have understood and was just waiting for Malcev to get it.
But then it turned out that there was a problem with speaker's argument.

During one of my visits, I attended a lecture by Andrey N.\ Kolmogorov himself.
During another visit, I talked to another great mathematician, Leonid V.\ Kantorovich.

\Qn Talked about what?

\Ar I had come across his approach to functional analysis and investigated some of the relevant algebraic notions \cite[\#10]{Gur}. In the process, I solved one of the open problems.
During my next visit to Novosibirsk, I made an appointment and 
went to meet Kantorovich. 
As I started to speak about the problem in question, he laughed and told me, as if confidentially: 
``I do not remember a thing about that.
Here's the address of Prof. Pinsker in Leningrad. Write to him.''
It turned out that the problem had been already solved by a Japanese mathematician.

Novosibirsk was also a great cultural center.
I heard there Bulat Okudzhava singing his songs and Bella Achmadulina reading her poetry.
And I made friends there, in particular Misha Taitslin.

\Qn What about Yuri Ershov?

\Ar We quickly developed very friendly relations, but then things changed dramatically.

\Qn Why?

\Ar I am not sure. He might have gone to the dark side.

\Qn That is confusing. Why do you think that you became friends in the first place?

\Ar As I mentioned, Ershov worked hard to understand my OAG stuff. I appreciated that.
When I came to Novosibirsk to arrange my candidate defense, NSU put me into a student dormitory.
Ershov had just left the dormitory, and I got his bed.
The two other guys in the room became our common friends.
The four of us played cards from time to time.

At a 1965 conference, Ershov protested when a speaker failed to mention my OAG classification among the important results involving logic and group theory. \looseness=-1
In 1968, he served as an opponent on my Doctor of Sciences defense \cite{W-DofS}.

\Qn And then you argued?

\Ar No. But various rumors started to reach me.
Apparently, the general atmosphere in Novosibirsk changed for the worse, especially in mathematics where antisemitism became rampant.

\Qn Did that touch you personally in any way?

\Ar Very little and indirectly.
One episode is related to ordered abelian groups.
After working for a few years on the so-called classical decision problem \cite{BGG}, I returned to OAGs.
In the algebraic literature, OAG theorems rarely were first-order. Typically they quantified over convex subgroups.
Accordingly, I investigated an expanded OAG theory where you can quantify over elements and over convex subgroups and which allowed me to formalize (in the expanded language) the OAG theorems in the algebraic literature.
Surprisingly, the expanded theory turned out to be easier to analyze than the first-order theory.
It remained decidable. Moreover, the decision algorithm became simpler. I submitted the resulting article to the Siberian Mathematical Journal, but the journal went mum for a long time.
My inquiries were unanswered. There were two relevant editors: Ershov and Kargapolov.
I wrote a private letter to Kargapolov and received a terse postcard reply ``This is not me."
Eventually I published the paper in the West \cite[\#25]{Gur}.

\Qn Hmm, this is all very circumstantial evidence.
You don't know what the referee or referees wrote.

\Ar I happened to know what one referee wrote.
Soon after the submission, I met Kokorin, and he volunteered that he wrote an enthusiastic referee report and that the paper should appear soon.
But you make a good point.

\Qn When you read good historians, as I like to do, something rubs off on you. Good historians like to see original documents.

\Ar There is a relevant document, a 1980 book Decidability Problems and Constructive Models, where Ershov uses the elementary theory of OAGs but avoids quoting my paper. ``In this section, we will study a number of questions related to elementary theories of OAGs.
A good effective classification of OAGs has not been found.
The results obtained in this section will be used in the following chapter to study the elementary theory of the $p$-adic field.'' \cite[Chapter~3, \S6]{Ershov}

\Qn My Russian is less than rudimentary.
It seems though that Ershov says something in the Introduction about authorship problems.

\Ar Yes, you are right. He ``decided not to distract the readers by discussing the problems (not always obvious) of authorship.''
He adds: ``The sources used by the author are given in the bibliography,'' but the 259-item bibliography has no sources on ordered abelian groups.

\Qn Did the book upset you?

\Ar A little, but not too much. I learned about the book in the middle of the Logic Year at the Hebrew University of Jerusalem (1980-81 academic year).
I was surrounded by the best model theorists and set theorists in the world.
I loved every minute of it.
The issues of decidability and constructive models seemed somewhat outdated.

\Qn How come?

\Ar Well, I was moving toward computer science.
In the 1960s, the gap between recursive-function decidability and feasibility (by which I mean practicality) did not attract much attention.
But the recursive-function decidability of a problem may not mean much in practice.
For example, Alfred Tarski produced a decision algorithm for the elementary theory of the field of real numbers \cite{Tarski}.
A feasible version of that algorithm would cause a revolution in program verification, but a feasible version is not known and probably does not exist.

Similarly, the undecidability/infeasibility gap did not attract attention for a while, but the recursive-function undecidability of a problem may not preclude practical solutions for practical problems.
For example, the halting problem is famously undecidable,
but the Terminator tool developed at Microsoft Research solves the halting (termination) problem for industrially relevant programs like Windows device drivers \cite{Cook,CPR}.

\Qn With all respect to practicality, let's put it aside and  concentrate on mathematics.
What's wrong, if anything, with recursive-function constructivity?

\Ar Recursive-function constructivity is just one kind of constructivity. But one size does not fit all, as they say.
Different parts of mathematics need different notions of constructivity.
I wrote about that in my critique of the recursive-function constructivism \cite[\#123]{Gur}.

\Qn Can you give me, here and now, an example of a different notion of constructivity?

\Ar G\"odel's constructible sets \cite{Goedel}.

\section{The back-and-forth method}

\Qn You said that Taimanov independently discovered the back-and-forth method in model theory.
Did somebody discover the method earlier?

\Ar Yes, the back-and-forth method was first discovered by Roland \Fraisse\ and published in his thesis \cite{Fraisse}.
The method gives a criterion for elementary equivalence.

I attended Taimanov's lecture where he presented the back-and-forth method. Even though his paper \cite{Taimanov} on that subject had been published, it might have been his first lecture on the back-and-forth method.
There were a lot of people and much interest in the lecture. I remember that Malcev was there.
In Taimanov's presentation, the method, which he called \ru{метод перекидки} and which may be translated as the \emph{shuttle method},  was rather dynamic.
And so was the presentation.
Usually reserved, Taimanov was excited and enthusiastic, and he moved energetically back and forth in front of the blackboard.
The presentation was rather technical and not easy to follow.

Mathematics often simplifies dynamic processes by making them static. Today the back-and-forth method can be described much more simply.

\Qn How?

\Ar Consider two relational structures $A$ and $B$ of the same signature. Define a decreasing sequence
\[ \s_0 \supseteq \s_1 \supseteq \s_2 \supseteq \dots \]
of sets $\s_n$ of isomorphisms $I$ from finite substructures of $A$ to finite substructures of $B$.
$\s_0$ is the set of all such finite isomorphisms.
Every $\s_{n+1}$ comprises the finite isomorphisms $I$ in $\s_n$ satisfying the following two conditions:

\smallskip
{\tt Forth: } For every element $a$ in $A$, there is a finite isomorphism $J\in\s_n$ that extends $I$ and whose domain contains $a$.

\smallskip
{\tt Back: } For every element $b$ in $B$, there is a finite isomorphism $J\in\s_n$ that extends $I$ and whose range contains $b$.

\smallskip\noindent
The structures $A$ and $B$ are called \emph{$n$-equivalent} if the empty isomorphism belongs to every $\s_n$.
It is not difficult to check that $A,B$ are elementarily equivalent if and only if they are $n$-equivalent for all $n$.

\smallskip
Apparently, the back-and-forth idea occurred to Taimanov after a visit to a textile factory and was suggested by the motion of a shuttle in a loom \cite{Taimanov}.

Andrzej Ehrenfeucht \cite{Ehrenfeucht2} gave a game-theoretic form to the original back-and-forth method.
Given $n$, consider an $n$-round game between Spoiler and Duplicator.
Intuitively, Duplicator wants to prove the claim that $A,B$ are $n$-equivalent while Spoiler wants to refute the claim.
To simplify the language, let's presume that Spoiler is male, and Duplicator female.

On the $k^{th}$ round, Spoiler first chooses one of the two structures. If $A$ is chosen then he chooses an element $a_k$ in $A$, and if $B$ is chosen then he chooses an element $b_k$ in $B$.
Then Duplicator chooses an element in the other structure.
If he chose an element $a_k$ in $A$, then she chooses an element $b_k$ in $B$, and if he chose an element $b_k$ in $B$, then she chooses an element $a_k$ in $A$.

After $n$ rounds, we usually have a mapping
\[ a_1\mapsto b_1, \dots, a_n\mapsto b_n; \]
otherwise (if $a_i=a_j$ but $b_i\ne b_j$) Spoiler already won.
If the mapping is an isomorphism from
the reduct $A\mathop{\upharpoonright}\{a_1,\dots,a_n\}$ of $A$ to
the reduct $B\mathop{\upharpoonright}\{b_1,\dots,b_n\}$ of $B$,
then Duplicator won; otherwise Spoiler won.
It is not difficult to prove that $A,B$ are $n$-equivalent if and only if Duplicator has a winning strategy in the $n$-round game.

\Qn Is there a good reference on the subject?

\Ar Yes, for example the book \cite{Hodges} by Wilfrid Hodges which is a standard reference; see \S3.2 there.
Nowadays few employ the original back-and-forth method.
People find the game-theoretic form easier to use.
One exception is the book \cite{Poizat} by Bruno Poizat.

\Qn Why do people find the game approach easier to use?

\Ar This is something psychological because mathematically the two approaches are equivalent.
Games provide a natural way to deal with quantifier alternation, even in simple cases.
Imagine that you need to explain to students that a function $f: \R\to \R$ approaches $b$ as $x$ approaches $a$. You say that for every real $\eps>0$ there exists a real $\delta>0$ such that for every real $x$, we have
\begin{equation}\label{ed}
0 < |x-a| < \delta \implies |f(x) - b| < \eps.
\end{equation}
A few mathematically gifted students understand, but most students have little idea what are you talking about.
What would you tell them?

\Qn I see where you are going. A natural thing to say is something like this:
\begin{equation}\label{game}
\begin{split}
&\text{If you give me some $\eps>0$,
I can come up with some $\delta>0$}\\
&\text{such that, whatever $x$ you choose, \eqref{ed} holds.}
\end{split}
\end{equation}
And \eqref{game} describes a game of sorts.

\Ar Furthermore, \eqref{game} asserts that you have a winning strategy in the game.
Strategies are complicated but intuitive objects.
It is easier for us to think of strategies than to deal with first-order claims with numerous quantifier alternations, like the claims that the empty isomorphism belongs to $\s_n$.
Notice that every step in the recursive definition of $\s_n$ introduces a quantifier alternation.

\section*{Acknowledgments}
Many thanks to Sem\"en Kutateladze, Stanislav Speranski, and especially Andreas Blass for useful critical remarks.

\end{document}